\documentclass{amsart}

\newtheorem{theorem}{Theorem}[section]
\newtheorem{lemma}[theorem]{Lemma}
\newtheorem{proposition}[theorem]{Proposition}

\newtheorem{conjecture}[theorem]{Conjecture}

\theoremstyle{definition}

\newtheorem{example}[theorem]{Example}

\theoremstyle{remark}
\newtheorem{remark}[theorem]{Remark}

\numberwithin{equation}{section}

\DeclareSymbolFont{AMSb}{U}{msb}{m}{n}
\DeclareMathSymbol{\F}{\mathbin}{AMSb}{"46}
\DeclareMathSymbol{\N}{\mathbin}{AMSb}{"4E}
\DeclareMathSymbol{\Z}{\mathbin}{AMSb}{"5A}
\DeclareMathSymbol{\R}{\mathbin}{AMSb}{"52}
\DeclareMathSymbol{\C}{\mathbin}{AMSb}{"43}

\begin{document} \title[Polynomial Recurrences and Cyclic Resultants]{Polynomial Recurrences and Cyclic Resultants}

\author{Christopher J. Hillar}
\address{Department of Mathematics, Texas A\&M University, College Station, TX 77843}
\email{chillar@math.tamu.edu}

\author{Lionel Levine}
\address{Department of Mathematics, University of California, Berkeley, CA
94720.}
\email{levine@math.berkeley.edu} \thanks{Both authors were supported under a NSF Graduate Research Fellowship.} 


% ----------------------------------------------------------------
\begin{abstract}
Let $K$ be an algebraically 
closed field of characteristic zero and let $f
\in K[x]$. The $m$-th {\it cyclic resultant} of $f$ is \[r_m =
\text{Res}(f,x^m-1).\]  A generic monic polynomial is 
determined
by its full sequence of cyclic resultants; however, the known
techniques proving this result give no effective computational
bounds.  We prove that a generic 
monic polynomial of degree $d$
is determined by its first $2^{d+1}$ cyclic resultants and that
a generic monic reciprocal polynomial of even degree $d$ is 
determined by its first $2\cdot 3^{d/2}$ of them.  In
addition, we show that cyclic resultants satisfy a polynomial
recurrence of 
length $d+1$.  This result gives evidence
supporting the conjecture of Sturmfels and Zworski that $d+1$
resultants determine $f$.  In the process, we establish two
general 
results of independent interest: we show that certain
Toeplitz determinants are sufficient to determine whether a
sequence is linearly recurrent, and we give conditions 
under
which a linearly recurrent sequence satisfies a polynomial
recurrence of shorter length.\end{abstract} \maketitle \section{Introduction}

Let $K$ be an algebraically 
closed field of characteristic zero.
Given a monic polynomial
    $$ f(x) = \prod_{i=1}^d (x - \lambda_i) \in K[x], $$
the $m$-th \textit{cyclic resultant} of $f$ is
    \begin{equation}
    \label{cyclicdefn}    r_m(f) = \text{Res}(f,x^m-1)
        = \prod_{i=1}^d (\lambda_i^m - 1).
    \end{equation} 

One motivation for the study of cyclic resultants comes from topological
dynamics.  Sequences of the form (\ref{cyclicdefn}) count periodic 
points for toral endomorphisms.  If $A = (
a_{ij})$ is a $d \times d$ 
integer matrix, then $A$ defines an endomorphism $T$ of the $d$-dimensional torus $\mathbb T^d = 
\R^d/\Z^d$ given by \[T(\mathbf{x}) = 
A\mathbf{x} \mod \Z^d.\]  Let 
$\text{Per}_m(T) = \{\mathbf{x} \in \mathbb T^{d} : T^m(\mathbf{x}) = 
\mathbf{x}\}$ be the set of points on the torus fixed by the map 
$T^m$.  Under the ergodicity condition that no eigenvalue of $A$ is a 
root of unity, it follows (see \cite{Ward}) that \[\#\text{Per}_m(T) = 
|\det(A^m-\text{Id})| = |r_m(f)|,\] in which $f$ is the characteristic 
polynomial of $A$.

In connection with number theory, cyclic resultants were also studied
by Pierce and Lehmer \cite{Ward} in the hope of using them to
produce large primes.  As a simple example, the Mersenne
numbers $M_m = 2^m - 1$ arise as cyclic resultants of the polynomial
$f(x) =x-2$.  Indeed, the map $T(x) = 2x \mod 1$ has precisely $M_m$
points of period $m$.  Further motivation comes from knot theory
\cite{knots}, Lagrangian mechanics \cite{guillemin,Zworski}, and,
more recently, in the study of amoebas of varieties \cite{purbhoo}
and quantum computing \cite{kedlaya}.

The problem of recovering a polynomial from its sequence of cyclic 
resultants arises naturally in several applications.  Commonly, an 
explicit bound $N=N(d)$ is desired in terms of the degree $d$ of $f$ so 
that the first $N$ resultants $r_1,\ldots,r_N$ determine $f$ (see \cite{Zworski, kedlaya}).  For instance, given a toral endomorphism of the 
type discussed above, one would like to use a minimal amount of (coarse) 
period data to recover the spectrum of the matrix $A$.  In general, 
reconstruction of a polynomial from its sequence of cyclic resultants 
seems to be a difficult problem.  While it is known \cite{chris} that in 
many instances the full sequence of resultants determines $f$, this result 
is of little use in computation.  One purpose of the present article is to 
give explicit upper bounds on the complexity of this problem.  Our main 
result in this direction is

\begin{theorem}\label{upperboundthm}
A generic monic polynomial $f(x) \in K[x]$ of degree $d$ is
determined by its first $2^{d+1}$ cyclic resultants $r_1, \dots,
r_{2^{d+1}}$.
A generic monic reciprocal polynomial of even degree $d$ is determined by its first $2\cdot 3^{d/2}$ cyclic resultants.
\end{theorem}

Emperical evidence suggests that Theorem \ref{upperboundthm} is
far from tight.  A conjecture of Sturmfels and Zworski addresses the
special case of a \emph{reciprocal} polynomial $f$, that is,
one satisfying $f(1/x) = x^d f(x)$.

\begin{conjecture}\label{sturmzworskiconj}
A reciprocal monic polynomial $f(x) \in K[x]$ of even degree $d$ is
determined by its first $d/2+1$ cyclic resultants.
\end{conjecture}

Recently, there has been some progress on this conjecture for a special
class of reciprocal polynomials.  Kedlaya \cite{kedlaya} has shown that
for a certain reciprocal polynomial $f$ of degree $d$ arising from 
the numerator $P(t)$ of a zeta function of a curve over a finite field 
$\F_q$, the first $d$ resultants are 
sufficient to recover $f$.  He uses this result to give a quantum 
algorithm that computes $P(t)$ in time polynomial in the degree of the 
curve and $\log q$.  A proof of Conjecture~\ref{sturmzworskiconj} would 
further reduce the running time for Kedlaya's algorithm.
We offer the following related conjecture.

\begin{conjecture}\label{modifiedconj}
A generic monic polynomial $f(x)
\in K[x]$ of degree $d$ is determined by its first $d+1$ cyclic
resultants. \end{conjecture}

Presently Conjecture \ref{modifiedconj} is verified only up to $d=4$ (see
Section \ref{computexperiments}); however, we are able to offer a result
in the direction of Conjectures \ref{sturmzworskiconj} and
\ref{modifiedconj}.  We say that a sequence $\{a_n\}_{n \geq 1}$, $a_n \in K$ obeys a {\it polynomial recurrence of length $\ell$} if there is a
polynomial $P \in K[x_1, \dots, x_\ell]$ such that $P(a_n, \dots,
a_{n+\ell-1}) = 0$ for all $n \geq 1$.  Our theorem may then be stated as
follows. 

\begin{theorem}\label{cyclicpolyrecthm} Let
$f \in K[x]$ be a monic polynomial of degree $d$.  The sequence
$\{r_n\}_{n \geq 1}$ of cyclic resultants of $f$ obeys a polynomial 
recurrence of length $d+1$.  Moreover if $f$ is assumed reciprocal of even degree $d$, then 
$\{r_n\}$ obeys a polynomial recurrence of length $d/2+1$. 
\end{theorem}

Explicit polynomial recurrences witnessing 
Theorem~\ref{cyclicpolyrecthm} in the cases $d=1,2$ may be found 
in Section~$5$.

The main tools used in our analysis are two general results relating
linear and polynomial recurrences.   Let $\{a_n\}_{n \geq 1}$ be given by
    \begin{equation}
    \label{generallinrec}
    a_n = \sum_{i=1}^{t} p_i(n) \mu_i^n,
    \end{equation}
where the $p_i$ are nonzero polynomials in $K[x]$, and the $\mu_i$ are nonzero elements of $K$.  Let $\ell_i = \text{deg}(p_i)+1$.  It is well-known (see, e.g.\ \cite{diffeq}) that the sequence $a$ obeys a linear recurrence of length $\ell = \ell_1 + \cdots + \ell_t$; namely,
	\[ a_{n+\ell} + c_1 a_{n+\ell-1} + \cdots + c_\ell a_n = 0,  \qquad n\geq 1 \]
with coefficients $c_i$ determined by
	\[ x^\ell + c_1 x^{\ell-1} + \cdots + c_n = \prod(x-\mu_i)^{\ell_i}. \]
We say that $a$ obeys a simple linear recurrence if all $\ell_i=1$.
 
The following result gives conditions under which a linearly recurrent sequence satisfies a polynomial recurrence of shorter length.

\begin{theorem}\label{linrecurthm}
Let $\{a_n\}_{n \geq 1}$ be given by {\em(\ref{generallinrec})}, and let $r$ be the rank of the multiplicative group $A \subset K^*$ generated by the $\mu_i$.  Then the sequence $\{a_n\}$ obeys a polynomial recurrence of 
length $r+2$.  Moreover, if $\{a_n\}$ satisfies a simple linear 
recurrence, then $\{a_n\}$ obeys a polynomial recurrence of length $r+1$.
\end{theorem}

\begin{example}
Consider the Fibonacci sequence $F_n = \frac{1}{\sqrt{5}} (\mu_+^n -
\mu_-^n)$, where $\mu_\pm = \frac{1 \pm \sqrt{5}}{2}$.  Since
$\mu_- = -\mu_+^{-1}$, the group $A$ in Theorem~\ref{linrecurthm} is
generated by $\mu_+$ and $-1$, hence it has rank $r=1$.  As $F_n$ obeys a simple linear recurrence, we expect the
sequence $\{F_n\}$ to obey a polynomial recurrence of length $r+1=2$.
Indeed, it is well-known and easily seen by induction that $F_n^2 - F_n F_{n-1} - F_{n-1}^2 = (-1)^n$, so every pair $(F_{n-1}, F_n)$ lies on the zero-locus of $P(x,y) = (y^2 - xy - x^2)^2 - 1$.
\end{example}

\begin{example}
Let $a_n = q(n)$ be a quasi-polynomial of degree $d$ and period $N$;
that is, there are polynomials $q_0, q_1, \ldots, q_{N-1} \in K[x]$ of
degree at most $d$, such that $q(n) = q_i(n)$ whenever $n \equiv i$ (mod
$N$).  It is elementary that any such sequence can be expressed in the
form (\ref{generallinrec}) with $t=d$ and $\mu_i = \zeta^{i-1}$, where
$\zeta$ is a primitive $d$-th root of unity.  Since the group $A \simeq
\Z/d\Z$ has rank zero, Theorem~\ref{linrecurthm} asserts that the
sequence $a_n$ obeys a polynomial recurrence of length $2$.  By contrast,
for suitably chosen $q$ the shortest linear recurrence for $a_n$ has
length $(d+1)N + 1$.
\end{example}

We remark that Theorem~\ref{linrecurthm} does not always give the
shortest length of a polynomial recurrence.  For instance, the sequence 
$ a_n = 2^n + (-2)^n + 3^n + (-3)^n$ satisfies the length-2 polynomial 
recurrence $a_n a_{n+1} = 0$,
while the theorem only guarantees the existence of a recurrence of
length $r+1 = 3$.

The second general result we use gives a polynomial recurrence of length
$2\ell-1$ which ``detects'' for the existence of a (homogeneous) linear
recurrence of length at most $\ell$.  Given a sequence $\{a_n\}_{n \geq
1}$, consider the $\ell \times \ell$ Toeplitz matrix
    \begin{equation}
    \label{Toeplitzdefn}
     A_{\ell,n} = \left[ \begin{array}{cccc}
    a_n & a_{n+1} & \ldots & a_{n+\ell-1} \\
    a_{n-1} & a_n & \ldots & a_{n+\ell-2} \\
    \vdots & \vdots & \ddots & \vdots \\
    a_{n-\ell+1} & a_{n-\ell+2} & \ldots & a_n \end{array} \right].
    \end{equation}

\begin{theorem}
\label{Toeplitz}
The sequence $\{a_n\}_{n \geq 1}$ satisfies a homogeneous linear
recurrence of length at most $\ell$ if and only if every Toeplitz
determinant $\det A_{\ell,n}$ vanishes, $n \geq \ell$.
\end{theorem}

Although this result appears to be known in some form (for example, it is implicit in the treatment of ``number walls'' in \cite{cg}) we include a proof in Section~2 as we were unable to find a reference.

All linear recurrence relations in this paper will henceforth be assumed
homogeneous.  

\begin{example}
In the case $\ell = 2$, the theorem asserts that $a_n$ is an exponential 
sequence $c \mu^n$ if and only if $a_n^2 = a_{n-1} a_{n+1}$ for all $n$.
\end{example} 

The rest of the paper is organized as follows.  In Section~2, we
prove Theorems~\ref{linrecurthm} and~\ref{Toeplitz}.  The proof the
former theorem reduces essentially to the computation of the Krull
dimension of a semigroup algebra, and that of the latter theorem to
an inductive application of Dodgson's rule.  In Section~3, we
establish a Toeplitz determinant factorization which will be used in
the proof of Theorem~\ref{upperboundthm}, along with some related
factorizations of independent interest.  In Section~4 we apply these
results, together with those of \cite{chris}, to prove Theorems~\ref{upperboundthm} and~\ref{cyclicpolyrecthm}.  Finally, in
Section~5, we present computational evidence supporting 
Conjecture~\ref{modifiedconj}.

We thank Bernd Sturmfels and Maciej Zworski for bringing this 
problem to our attention and for useful discussions.

\section{Linear and Polynomial Recurrences}

Let $S$ denote the collection of all sequences $\{a_n\}_{n \geq 1}$ with terms in $K$.  Pointwise sum and product give $S$ the structure of a commutative $K$-algebra with unit.  We denote by $E : S \rightarrow S$ the $K$-algebra endomorphism $(Ea)_n = a_{n+1}$ (the ``shift operator'').  

For $\xi \in K^*$ denote by $e(\xi)$ the exponential sequence $e(\xi)_n = \xi^n$; note that $e(1)$ is the unit element of $S$.  We will make use of the fact that for distinct $\xi_1, \ldots, \xi_m$ the sequences $e(\xi_i)$ are linearly independent over $K$ (the determinant $|e(\xi_i)_j|_{i,j=1}^m$ is Vandermonde).  Denote by $\delta$ the sequence $\delta_n = n$.  Then a sequence of the form (\ref{generallinrec}) can be expressed 
	\begin{equation}
	\label{restatementwithoutn}
	a = \sum_{i=1}^t p_i(\delta) e(\mu_i).
	\end{equation}
	
The proof of Theorem \ref{linrecurthm} will make use of the subalgebra \[R = K [a, Ea, E^2 a, \ldots ] \subset S\] generated by the sequence $a$ together with its leftward shifts.  This is a finitely generated $K$-algebra because $a$ obeys a linear recurrence.

Let $Q$ be a commutative semigroup, written multiplicatively.  The
{\it semigroup algebra} $K[Q]$ has $K$-basis indexed by the elements
of $Q$.  The basis element corresponding to $\mu \in Q$ is written 
$[\mu]$.  Multiplication is defined on basis elements by $[\mu] [\lambda] 
= [\mu \lambda]$ and extended by linearity.

\begin{lemma} \label{isom}
Let $a_n$ be given by {\em (\ref{generallinrec})},  and let $Q \subseteq K^*$  be the multiplicative semigroup generated by $\mu_1,\ldots,\mu_t$.  

\begin{enumerate}
\item There is an inclusion of $K$-algebras $R \hookrightarrow K[Q][x]$.  
\item If the sequence $a_n$ obeys a simple linear recurrence, there is an isomorphism of $K$-algebras $R \simeq K[Q]$. 
\end{enumerate}
\end{lemma}

\begin{proof} 
Write $\ell = \sum (\text{deg} ~ p_i + 1)$.  Since $a$ obeys a linear recurrence of length $\ell+1$, we have    
	\begin{equation}     \label{thebasics}   
	R = K [a, Ea, \ldots, E^{\ell-1} a].    
	\end{equation} 
Since $E e(\xi) = \xi e(\xi)$ and $E\delta = \delta + 1$, by (\ref{restatementwithoutn}) we have   
	\begin{equation} \label{alltermsatonce}     
	E^j a = \sum_{i=1}^t p_i(\delta+j) \mu_i^j e(\mu_i),    
	\end{equation} 
and hence there is an inclusion $R \subseteq R'[\delta]$, where
	\begin{equation}    \label{almostthere}   
	R' = K [e(\mu_1), \ldots, e(\mu_t)]. 
	\end{equation}
Since $e(\mu_i) e(\mu_j) = e(\mu_i \mu_j)$ the linear map $e: K[Q] \rightarrow R'$ sending $[\mu_i] \mapsto e(\mu_i)$ is a $K$-algebra homomorphism.  Since the exponential sequences $e(\mu_i)$ are linearly independent, it is an isomorphism.  Thus $R' \simeq K[Q]$.

If $\{a_n\}$ satisfies a simple linear recurrence, then the polynomials $p_i(n)$ are nonzero constants.  From (\ref{alltermsatonce}) we have
	$$ E^j a = \sum_{i=1}^t p_i \mu_i^j e(\mu_i). $$
Thus the linear span of the sequences $a$, $Ea,\ldots,E^{\ell-1} a$ coincides with that of $e(\mu_1), \ldots, e(\mu_\ell)$ (the transition matrix is the product 
of a Vandermonde and an invertible diagonal matrix).  It follows from (\ref{thebasics}) that
	\begin{equation} 
	\label{expos}   R = K [e(\mu_1), \ldots, e(\mu_t)] = R' \simeq K[Q]. 
	\end{equation}
This completes the proof of (2).

To prove (1), it suffices to show that $\delta$ is transcendental over $R'$.  
Suppose it had algebraic degree $m$.  From among the algebraic relations
	\begin{equation}\label{linrelationinduction}
	\rho_0 \delta^m + \rho_1 \delta^{m-1} + \cdots + \rho_m = 0,       \qquad \rho_j \in R', ~\rho_0 \neq 0,
	\end{equation}
writing 
	$$\rho_0 = \sum_{i=1}^s b_i e(\xi_i), \qquad b_i \in K^*, ~\xi_i \in Q,$$ 
choose a relation in which the number of terms $s$ in the leading coefficient $\rho_0$ is minimal.  This minimality forces $e(\xi_1), \ldots, e(\xi_s)$ to be linearly independent; equivalently, the $\xi_i$ must be distinct.  Multiplying (\ref{linrelationinduction}) by $b_s^{-1} e(\xi_s^{-1})$ we may assume that $b_s=\xi_s=1$.
 
Suppose first that $s>1$.  Applying the difference operator $\Delta=E-1$ to (\ref{linrelationinduction}) 
we obtain
	$$ E (\rho_0) (\delta+1)^m - \rho_0 \delta^m + q(\delta) = 0, $$
where $q \in R'[x]$ is a polynomial of degree at most $m-1$.  Thus
	\begin{eqnarray*}
	0 &=& \sum_{i=1}^s b_i \xi_i e(\xi_i) (\delta+1)^m - \sum_{i=1}^s b_i e(\xi_i) \delta^m + q(\delta) \\
	  &=& \sum_{i=1}^{s-1} b_i (\xi_i - 1) e(\xi_i) \delta^m + q(\delta) + \tilde{q}(\delta),
	\end{eqnarray*}
in which $\tilde{q} \in R'[x]$ again has degree at most $m-1$.  By the minimality of $s$, the coefficient of $\delta^m$ must vanish, and this contradicts the linear independence of the sequences $e(\xi_i)$. 

It remains to consider the case $s=1$.  By our rescaling convention, $b_1 = \xi_1 = 1$, hence $\rho_0 = e(1)$.  Writing $\rho_1 = \sum_{i=1}^{s'} c_i e(\psi_i)$ for distinct $\psi_i$, applying $\Delta$ to the relation (\ref{linrelationinduction}), we obtain
	$$ e(1) \left[ (\delta+1)^m - \delta^m \right] + \sum_{i=1}^{s'} c_i \left[ \psi_i e(\psi_i) (\delta+1)^{m-1} - e(\psi_i) \delta^{m-1} \right] + p(\delta) = 0, $$
where $p \in R'[x]$ is a polynomial of degree at most $m-2$.  The coefficient of $\delta^{m-1}$ is
	\begin{equation} \label{seeya} 
	m e(1) + \sum_{i=1}^{s'} (\psi_i - 1) e(\psi_i) = 0. 
	\end{equation}
If $s'=0$ we obtain $me(1)=0$, a contradiction as $K$ has characteristic $0$; if $s'=1$ and $\psi_1=1$ we obtain the same contradiction.  Finally, if some $\psi_i\neq 1$, then (\ref{seeya}) contradicts the linear independence of the sequences $e(\psi_i)$ and $e(1)$.
\end{proof}
	
\begin{lemma}
\label{equalities}
Let 
$A$ be a finitely generated abelian group, and fix a set of generators
$q_1, \ldots, q_\ell$ for $A$.  Let $Q \subset A$ be the semigroup
generated by the $q_i$.  The 
following are equal.
\begin{enumerate}
\item The rank of $A$;
\item The Krull dimension of $K[Q]$;
\item The maximum number of elements in $K[Q]$ algebraically independent
over $K$.
\end{enumerate}
\end{lemma}

\begin{proof}
Write $A = B \oplus C$, with $B$ finite and $C$ free abelian.  By
Maschke's Theorem, $K[B]$ is a finite product of 
copies of $K$, hence
$K[A] \simeq K[B] \otimes_K K[C]$ is a finite product of copies of $K[C]$.
Since $K[C]$ is a Laurent polynomial ring, it follows that $K[A]$ is
reduced, and hence $K[Q]$ is reduced.  Now by \cite[p. 466]{Cox1}, the
maximum number of algebraically independent elements in $K[Q]$ is equal to
its Krull dimension.   
Finally, by Proposition 7.5 in \cite{CCA}, the Krull
dimension of $K[Q]$ is equal
to the rank of $A$.
\end{proof}

\begin{proof}[Proof of Theorem~\ref{linrecurthm}]
By 
Lemma~\ref{isom}, we have $R \hookrightarrow K[Q][x] \simeq K[Q \times
\N]$.  
Thus by Lemma~\ref{equalities}, the maximum number of algebraically
independent elements in 
$R$ is at most $\text{rank}(A \times \Z) = r+1$.
In particular, the $r+2$ elements $a, Ea, \ldots, E^{r+1}a \in R$ are
algebraically dependent over $K$; that is, the 
sequence $a$ obeys a
polynomial recurrence of length $r+2$.

If $a$ satisfies a simple linear recurrence, then $R \simeq K[Q]$ by
Lemma~\ref{isom}, and the elements $a, Ea, 
\ldots, E^r a$ are
algebraically dependent over $K$, so $a$ obeys a polynomial recurrence of
length $r+1$.
\end{proof}

We now turn to the proof of Theorem~\ref{Toeplitz}.  
The key step uses
Dodgson's rule \cite{Zeil} relating the determinant of an $\ell \times
\ell$ matrix to its four corner $\ell-1 \times \ell-1$ minors and its
central 
$\ell-2 \times \ell-2$ minor.  For Toeplitz matrices the rule
assumes a particularly simple form: with $A_{\ell,m}$ defined as in
(\ref{Toeplitzdefn}), we have for $2 \leq \ell \leq m$
\begin{equation}    \label{Dodgson}    \det A_{\ell,m} \det A_{\ell-2, m}
    = (\det A_{\ell-1, m})^2 - \det A_{\ell-1, m-1} \det A_{\ell-1, m+1}.    \end{equation}
\noindent {\it Proof of Theorem~\ref{Toeplitz}.}
We induct on $\ell$.  If $x_1 a_n + x_2 a_{n+1} + \dots + x_\ell
a_{n+\ell-1} = 0$ for all $n \geq 1$, then the vector $[x_1, \dots,
x_\ell]^T$ lies in the kernel of every $A_{\ell,n}$.  Conversely,
suppose $A_{\ell,n}$ is singular for all $n \geq \ell$.  For each such
$n$, let $\textbf{x}_n = [x_{n1}, \dots, x_{n\ell}]^T$ be a nonzero
vector in the kernel of $A_{\ell,n}$, and consider the $2 \times
\ell$ matrices \[X_n = \left[\begin{array}{c}\textbf{x}_n^{T} \\
\textbf{x}_{n+1}^{T} \end{array}\right]=
\left[\begin{array}{ccc}x_{n1} & \cdots & x_{n\ell} \\x_{n+1,1} &
\cdots & x_{n+1,\ell}\end{array}\right]. \] If every $X_n$ has
rank one, then each vector $\textbf{x}_n$ is a scalar multiple of
$\textbf{x}_\ell$, and $\{a_n\}$ satisfies the linear recurrence
$x_{\ell1} a_n + \dots + x_{\ell\ell} a_{n+\ell-1} = 0$.

Suppose now that some $X_n$ has rank two. The transposes of the first
$\ell-1$ row vectors of $A_{\ell,n}$ all lie in the kernel of $X_n$,
so they must be linearly dependent.  In particular, the upper-left
minor $\det A_{\ell-1, n}$ vanishes.  We now induct forwards and
backwards on $m$ to show that $\det A_{\ell-1, m}$ vanishes for all
$m \geq \ell$.  The left hand side of (\ref{Dodgson}) is zero since
$A_{\ell,m}$ is singular.  Therefore, if either $\det A_{\ell-1,
m-1}$ or $\det A_{\ell-1, m+1}$ vanishes, then $\det A_{\ell-1, m}$
must vanish as well.  This completes the induction on $m$.  By
induction on $\ell$, the sequence $\{a_n\}_{n \geq 2}$ satisfies a
linear recurrence of length at most $\ell-1$.  This trivially implies
that $\{a_n\}_{n \geq 1}$ satisfies a linear recurrence of length at
most $\ell$. \qed

\section{Determinant Factorizations}

If $\{a_n\}_{n \geq 1}$ satisfies a simple linear recurrence of
length $\ell+1$, the determinants $\det A_{\ell,n}$ have a simple
closed form.  Note that $\ell+1$ is the minimum length for which
these determinants do not vanish, by Theorem~\ref{Toeplitz}.

\begin{theorem}
\label{determinant}
Let $\mu_i, c_i \in K^*$ {\em (}$i=1,\ldots,\ell${\em )}, and set
$a_n = \sum_{i=1}^\ell c_i \mu_i^n$.  The determinant of the $\ell
\times \ell$ Toeplitz matrix $A_{\ell,n}$ defined in
$(\ref{Toeplitzdefn})$ has the factorization
    \begin{equation}
    \label{factorization}
    \det A_{\ell,n} = (-1)^{\ell(\ell-1)/2} c_1 \cdots c_\ell
    (\mu_1 \cdots \mu_\ell)^{n-\ell+1}
    \prod_{i < j} (\mu_i-\mu_j)^2.
    \end{equation}
\end{theorem}

\begin{proof}
Consider the Vandermonde matrices $V=(\mu_i^{j-1})_{i,j=1}^\ell$ and
$V'=(\mu_j^{\ell-i})_{i,j=1}^\ell$, and let $D$ be the diagonal
matrix with diagonal entries $c_i \mu_i^{n-\ell+1}$, $i=1, \ldots,
\ell$.  The $(i,j)$-entry of the product $V'DV$ is then
    $$ (V'DV)_{i,j} = \sum_{k=1}^\ell \mu_k^{\ell-i} c_k
    \mu_k^{n-\ell+1} \mu_k^{j-1} = \sum_k c_k \mu_k^{n-i+j}
        = a_{n-i+j}. $$
Thus $A_{\ell,n} = V'DV$.  Since $V'$ differs from $V^T$ by a row
permutation of length $\ell(\ell-1)/2$ we obtain
    $$ \det A_{\ell,n} = (-1)^{\ell(\ell-1)/2} \det D (\det V)^2, $$
which yields (\ref{factorization}).
\end{proof}

It seems likely that Theorem~\ref{determinant} can be extended to
cover the situation of an arbitrary linear recurrence.  One
difficulty is that in general, the factors $\mu_i - \mu_j$ occur with
multiplicity.  For example, if \[a_n = (bn^3+cn^2+dn+e)\mu_1^n +
(fn+g)\mu_2^n + h\mu_3^n,\]  then \[\det A_{7,n} = 1296 b^4 f^2 h
\mu_1^{4n-12} \mu_2^{2n-8} \mu_3^{n-6} (\mu_1-\mu_2)^{16} (\mu_1-\mu_3)^8
(\mu_2-\mu_3)^4.\]  Of particular interest is the special case when $a_n$
is a polynomial function of $n$.

\begin{proposition}\label{d!prop}
Let $p(x)$ be a monic polynomial of degree $\ell$.  As a polynomial
in $x$, the $(\ell+1) \times (\ell+1)$ Toeplitz determinant
$\Delta(x) = \det(p(x-i+j))_{i,j=0}^\ell$ equals the constant
$\ell!^{\ell+1}$, independent of $p$.
\end{proposition}

The proof of this proposition will make use of the following
generalization of the Vandermonde determinant \cite{krat}.

\begin{lemma}\label{vandlem}
Let $p_j(x)$ be polynomials {\em(}$j=0,\ldots,\ell${\em)} with leading
coefficients $a_j$ and $\deg(p_j) = \ell-j$.  If $y_0,\ldots,y_\ell$
are indeterminates, then
\begin{equation*} 
\det \left(p_j(y_i) \right)_{i,j=0}^{\ell} \ = \ a_0 \cdots a_\ell
\prod_{i<j}(y_i-y_j).
\end{equation*}
\end{lemma}

\begin{proof}
Perform elementary column operations from right to left,
reducing the determinantal calculation to that of the Vandermonde
identity.
\end{proof}

\begin{proof}[Proof of Proposition \ref{d!prop}]
Set $q_1(x) = \frac{p(x+1) - p(x)}{\ell}$, which is monic of
degree $\ell-1$. Subtracting the first column from the second, the
second from the third, and 
so on, we obtain \[\Delta(x) = \ell^\ell
\left|\begin{array}{cccc}p(x) & q_1(x) & \cdots
& q_1(x+\ell-1) \\p(x-1) & q_1(x-1) & \cdots & q_1(x+\ell-2) \\\vdots
& \vdots & 
\ddots & \vdots \\p(x-\ell) & q_1(x-\ell) & \cdots &
q_1(x-1)\end{array}\right|.\]

Next, set $q_{i+1}(x) = \frac{q_{i}(x+1)-q_{i}(x)}{\ell-i}$ for $i =
1,\ldots, 
\ell-1$ and repeat the above reductions, initiating the
column operations at the $(i+1)$-st column.  At the conclusion of this
process, we end up with the determinant 
\[\Delta(x) = \ell^\ell
(\ell-1)^{\ell-1} \cdots 1^1 \left|\begin{array}{cccc}p(x) & q_1(x) &
\cdots & q_\ell(x) \\ p(x-1) & q_1(x-1) & \cdots & q_\ell(x-1)
\\ \vdots & 
\vdots & \ddots & \vdots \\ p(x-\ell) & q_1(x-\ell)  &
\cdots & q_\ell(x-\ell) \end{array} \right|. \]  Finally, setting $y_i =
x-i$ in the statement of  Lemma 
\ref{vandlem}, it follows that  \[\Delta(x)
= \prod_{k=1}^\ell k^k
\cdot \prod_{0 \leq i < j \leq \ell} (j-i) = \prod_{k=1}^\ell k^k
\cdot \prod_{k=1}^\ell k^{\ell+1-k} = 
\ell!^{\ell+1}. \qed \]
\renewcommand{\qedsymbol}{}
\end{proof}

Applying Theorem \ref{determinant} to the cyclic resultants $r_n$ given in
(\ref{cyclicdefn}) we obtain the 
following.

\begin{proposition}
\label{silly}
The determinant of the $2^d \times 2^d$ Toeplitz matrix $R_n =
(r_{n-i+j})_{i,j=1}^{2^d}$ has the factorization
    
\begin{equation}\label{rntoeplitz}
    \det R_n = (\lambda_1 \cdots \lambda_d)^{2^{d-1} (n - 2^{d-1})}
    \prod_{S,T} (\lambda_S - \lambda_T)^{2^{d-|S \cup T|}},
\end{equation}
where $\lambda_S := \prod_{i \in S} \lambda_i$, and the product is taken 
over all ordered pairs of disjoint subsets $S,T \subseteq \{1, \dots, d\}$,
not both empty. 
\end{proposition}

\begin{proof} Expanding the product (\ref{cyclicdefn}) yields 
$r_n = \sum_{S \subseteq \{1, \ldots, d\}} (-1)^{d-|S|}
\lambda_S^n$.  Since $\prod_S \lambda_S = (\lambda_1 \cdots
\lambda_d)^{2^{d-1}}$, by Theorem \ref{determinant} with $\ell =
2^d$ we have
    \begin{equation}
    \label{alternateform}
    \det R_n = (\lambda_1 \cdots
        \lambda_d)^{2^{d-1} (n-2^d+1)}
        \prod_{S \neq T} (\lambda_S - \lambda_T).
    \end{equation}
For each pair of subsets $S,T \subseteq \{1, \ldots, d\}$ write $S' = S - S
\cap T$, $T' = T - S \cap T$.  Then $\lambda_S - \lambda_T = \lambda_{S
\cap T} (\lambda_{S'} - \lambda_{T'})$.  Each ordered pair $(S',T')$ of
disjoint subsets of $\{1, \dots, d\}$, not both empty, arises in this
way from $2^{d-|S' \cup T'|}$ different ordered pairs $(S,T)$, $S
\neq T$.  The product in (\ref{alternateform}) can therefore be
rewritten as \[\prod_{S \neq T} (\lambda_S - \lambda_T) = \prod_{S
\neq T}  \lambda_{S \cap T} \prod_{ \begin{array}{c} {\scriptstyle S'\cap T' = \emptyset} \\ {\scriptstyle (S',T') \neq (\emptyset,\emptyset) }
\end{array} }  (\lambda_{S'} - \lambda_{T'})^{2^{d-|S' \cup T'|}}.\]
Every element $i \in \{1, \dots, d\}$ is contained in $4^{d-1}-2^{d-1}$
intersections of the form $S \cap T$, $S \neq T$.  Thus, the exponent
on $\lambda_1 \cdots \lambda_d$ is $2^{d-1} (n - 2^d + 1) + 4^{d-1} -
2^{d-1} = 2^{d-1} (n - 2^{d-1})$.
\end{proof}

\begin{remark}
The degree of $\det R_n$ as a polynomial in the roots $\lambda_1,
\dots, \lambda_d$ is significantly smaller than might be expected.  
Every term in the expansion of $R_{2^d}$ has degree $d 4^d$.  On the other
hand, by (\ref{alternateform}) the degree of $\det R_{2^d}$ is only
$M(d) := \sum_{S,T \subseteq 
\{1, \dots, d\}} \text{max}\{|S|,|T|\}$.
It is easily seen (e.g. by the weak law of large numbers) that $M(d)
\sim \frac12 d 4^d$ as $d \rightarrow \infty$.  For large $d$, so
much cancellation occurs in the expansion of $\det R_{2^d}$
that its degree is approximately cut in half.
\end{remark}

\section{Cyclic Resultants}

We are now in a position to prove the first two results stated in
the introduction.  Expanding the product formula (\ref{cyclicdefn}) yields
    \begin{equation} \label{expansion}
    r_n = \sum_{S \subset \{1, \ldots, d\}} (-1)^{d-|S|}\lambda_S^n
    \end{equation}
where $\lambda_S = \prod_{i \in S} \lambda_i$.  Thus the sequence $r_n$
obeys a simple linear recurrence of length $2^d$.  Note, however, that the
coefficients are functions of the $\lambda_i$.

If $f$ is reciprocal of even degree $d$, its roots come in reciprocal pairs $\lambda_1^{\pm 1}, \ldots, \lambda_{d/2}^{\pm 1}$, and the expansion of (\ref{cyclicdefn}) yields
	\begin{equation*} \label{reciprocalexpansion}
	r_n = \sum_{S,T} 2^{\frac{d}{2}-|S\cup T|} (-1)^{|S\cup T|} \left(\frac{\lambda_S}{\lambda_T} \right)^n, \end{equation*}
where the sum is over disjoint subsets $S,T \subset \{1,\ldots,d/2\}$.

\begin{proof}[Proof of Theorem \ref{upperboundthm}]
We take ``generic'' to mean that $f$ does not have a root of unity as a 
zero and that no two products of distinct subsets of 
the roots of $f$ are equal; that is, $\lambda_S \neq \lambda_T$ for 
distinct subsets $S, T \subset \{1, \ldots, d\}$.  The first author, in Corollary~1.7 of 
\cite{chris}, extending work of Fried \cite{fried}, proves that a generic monic 
polynomial $f(x) \in K[x]$ is determined by its full 
sequence $\{r_n\}_{n \geq 1}$ of cyclic resultants.  By (\ref{expansion}) 
and Theorem~\ref{Toeplitz}, the sequence $\{r_n\}$ obeys a polynomial 
recurrence of length $2^{d+1}+1$ given by the vanishing of $(2^d + 1) 
\times (2^d + 1)$ Toeplitz determinants.  Unlike the linear recurrence of 
length $2^d + 1$, this recurrence is {\it independent} of the polynomial 
$f$.  By minor expansion along the bottom row, the recurrence may be 
expressed in the form
    \begin{equation} \label{rnRneqn}
    (\det R_n) \cdot r_{n+2^d+1}
        = P(r_{n-2^d+1}, \dots, r_{n+2^d}), \qquad n \geq 2^d
    \end{equation}
for some polynomial $P \in K[x_1, \ldots, x_{2^{d+1}}]$, where $R_n =
(r_{n+i-j})_{i,j=1}^{2^d}$.

By Proposition~\ref{silly}, a generic polynomial $f$ gives rise
to nonsingular $R_n$ for all $n$.  Moreover, a straightforward
calculation using equation (\ref{rntoeplitz}) reveals that \[\det R_n
= \frac{ (\det R_{2^d+1})^{n-2^d} }{ (\det R_{2^d})^{n-2^d-1}} .\]
It follows from this and (\ref{rnRneqn}) that for $m
\geq 2^{d} + 2^{d} + 1 = 2^{d+1}+1$, the resultant $r_m$ is
determined by the resultants $r_i$ with $i < m$.  In particular, the
values $r_1, \dots, r_{2^{d+1}}$ determine the full sequence of resultants, and hence they determine $f$.

If $f$ is reciprocal of even degree $d$, we take ``generic'' to mean that 
$f$ does not have a root of unity as a zero and that no two 
quotients of the form $\lambda_S/\lambda_T$ are equal, where $S$ and 
$T$ are disjoint subsets of $\{1,\ldots,d/2\}$.  As there are $3^{d/2}$ such pairs $(S,T)$, by 
Theorem~\ref{Toeplitz}, the sequence $\{r_n\}$ obeys a polynomial 
recurrence of length $2\cdot 3^{d/2}+1$ given by the vanishing of $(3^{d/2} + 1) 
\times (3^{d/2}+ 1)$ Toeplitz determinants.  Moreover, the matrices $R'_n = (r_{n-i+j})_{i,j=1}^{3^{d/2}}$ 
are nonsingular by Theorem~\ref{determinant}.  In \cite[Corollary 1.12]{chris}, it is shown that 
the sequence of nonzero cyclic resultants generated by a reciprocal polynomial of degree $d$ 
determines it. The proof now proceeds by minor expansion as before.
\end{proof}

\begin{proof}[Proof of Theorem \ref{cyclicpolyrecthm}]
The elements $\lambda_S \in K$ appearing in (\ref{expansion}) lie in the
multiplicative subgroup $A \subseteq K^*$ generated by $\lambda_1,
\ldots,\lambda_d$.  By Theorem~\ref{linrecurthm} applied to
(\ref{expansion}), the sequence $\{r_n\}$ obeys a polynomial recurrence of
length $\text{rank}~A + 1 \leq d+1$.

If $f(x)$ is reciprocal, its roots $\lambda_i$ come in reciprocal pairs,
and the rank of $A$ is at most $d/2$.  By
Theorem~\ref{linrecurthm} it follows that $\{r_n\}$ obeys a polynomial
recurrence of length $d/2+1$.
\end{proof}

We close this section with an example illustrating the techniques used in
proving Theorem \ref{upperboundthm}.  By (\ref{expansion}) the resultants
$r_m$ of a monic quadratic polynomial $f(x) = x^2 + ax + b$ satisfy a
linear recurrence of length~$5$.  By Theorem~\ref{Toeplitz} it follows
that \[ \det A_{5,n} = \left|\begin{array}{ccccc}r_n & r_{n-1} & r_{n-2} &
r_{n-3} & r_{n-4} \\ r_{n+1} & r_n & r_{n-1} & r_{n-2} & r_{n-3} \\r_{n+2}
& r_{n+1} & r_n & r_{n-1} & r_{n-2} \\ r_{n+3} & r_{n+2} & r_{n+1} & r_n &
r_{n-1} \\ r_{n+4} & r_{n+3} & r_{n+2} & r_{n+1} & r_n \end{array}\right|
= 0.\] Minor expansion along the bottom row gives our polynomial
recurrence
\begin{equation}\label{examplerecur}
r_{n+4} \left| \begin{array}{cccc} r_{n-1} & r_{n-2} & r_{n-3} & r_{n-4}
\\ r_n & r_{n-1} & r_{n-2} & r_{n-3} \\ r_{n+1} & r_n & r_{n-1} & r_{n-2}
\\ r_{n+2} & r_{n+1} & r_n & r_{n-1} \end{array} \right| = P(r_{n-4},
r_{n-3}, \ldots, r_{n+3}) \end{equation}
expressing $r_{n+4}$ as a rational function in $r_{n-4}, \ldots,
r_{n+3}$.  For generic $f$, the determinant on the left side of
(\ref{examplerecur}) never vanishes, so the entire sequence of resultants
is determined by the values of $r_1, \ldots, r_8$.

While (\ref{examplerecur}) gives a polynomial recurrence for $r_m$ of
length~$9$, there in fact exists a recurrence of length~$3$ by
Theorem~\ref{cyclicpolyrecthm}.  Explicitly, this recurrence is given in 
the next section.  The coefficients of this shorter 
recurrence depend on the polynomial $f$, whereas (\ref{examplerecur}) 
gives a universal recurrence independent of $f$. 

\section{Computations}\label{computexperiments}
We list here explicit polynomial recurrences witnessing 
Theorem~\ref{cyclicpolyrecthm} in the cases $d=1$ and $d=2$.  The cyclic 
resultants of a monic linear polynomial 
$f(x) = x + a$ obey the 
length-$2$ recurrence
	$$ r_{n+1} = -a r_n - a - 1. $$
For a monic quadratic $f(x) = x^2 + ax + b$ we have the length-$3$ 
polynomial recurrence
\begin{equation*}\label{polyrecd2}
\begin{split}
& (a+b+1) \left[ (a-2)r_{n+2} + a(a-b-1)r_{n+1} + (a-2b)b r_n -
(a-b-1)(a+b+1) \right] \\
& = -r_{n+2}^2 - 
(a-2b)r_{n+1}r_{n+2} + ab r_n r_{n+2} + (a-b-1)b r_{n+1}^2   \\
& \ \ \ \ - (a-2)b^2 r_n r_{n+1} - b^3 r_n^2.\\\end{split}
\end{equation*}

We close with the explicit 
reconstructions of polynomials from their 
cyclic resultants in the cases $d=2,3$.   For quadratic $f =  x^2 + a x + b$, two nonzero resultants 
suffice to determine $a$, $b$: 
\begin{equation*}
\begin{split} a = \frac{r_1^2-r_2}{2r_1}, \ ~b = \frac{r_1^2-2r_1+r_2}{2r_1}.  \\ \end{split} \end{equation*} For cubic $f =  x^3 + a x^2 + 
b x + c$, four nonzero resultants give inversion: \begin{equation*} \begin{split} a = &\ 
\frac{-12r_2r_1^3-12r_1r_2^2+3r_2^3-r_2r_1^4-8r_2r_1r_3+6r_1^2r_4}{24r_2r_1^
2},  b =  \frac{-r_1^2-2r_1+r_2}{2r_1}, \\  c = & \ 
\frac{-3r_2^3+r_2r_1^4+8r_2r_1r_3-6r_1^2r_4}{24r_1^2r_2}. \\
\end{split}
\end{equation*} In addition, a monic quartic polynomial can be explicitly reconstructed 
using five 
resultants; however, the expressions are too cumbersome to 
list here.  We remark that the pattern of monomial denominators found in the inversions above does not continue for higher degree reconstructions.

% ----------------------------------------------------------------


\begin{thebibliography}{99}

\bibitem{cg} J. H. Conway and R. Guy, \emph{The Book of Numbers}, Springer-Verlag, 1996.

\bibitem{Cox1}
D. Cox, J. Little, and D. O'Shea, \emph{Ideals, varieties, and
  algorithms}, second ed., Undergraduate Texts in Mathematics,
  Springer--Verlag, New York, 1997.

\bibitem{duistermaat}
J.J. Duistermaat and V. Guillemin, \emph{The spectrum of positive
elliptic operators and periodic bicharacteristics}, Inv. Math. 25
(1975) 39-79.

\bibitem{diffeq} S. Elaydi, \emph{An Introduction to Difference Equations}, Springer, 1999.

\bibitem{Ward}
G. Everest and T. Ward. Heights of Polynomials and Entropy in
Algebraic Dynamics. Springer-Verlag London Ltd., London, 1999.

\bibitem{fried}
D. Fried, \emph{Cyclic resultants of reciprocal polynomials}, in
Holomorphic Dynamics (Mexico 1986), Lecture Notes in Math. 1345,
Springer Verlag, 1988, 124-128.

\bibitem{guillemin}
V. Guillemin, \emph{Wave trace invariants}, Duke Math. J. 83
(1996), 287-352.

\bibitem{chris} C. Hillar, \emph{Cyclic resultants}, J. Symb. Comp. {\bf 39} (2005), 653--669; 
erratum, ibid. {\bf 40} (2005), 1126Ð-1127.

\bibitem{Zworski}
A. Iantchenko, J. Sj\"{o}strand, and M. Zworski, \emph{Birkhoff
normal forms in semi-classical inverse problems}, Math. Res. Lett. 9 (2002), 
337-362.

\bibitem{kedlaya}
K. Kedlaya, \emph{Quantum computation of zeta functions of curves}, 
Computational Complexity \textbf{15} (2006), 1--19.  MR2226067

\bibitem{krat}
C. Krattenthaler, ``Advanced Determinant Calculus,'' \emph{Sem. Lothar. Combin.} 
{\bf 42} (1999), Art. B42q. 

\bibitem{CCA}
E. Miller and B. Sturmfels, \emph{Combinatorial Commutative Algebra}, Springer, 2004.

\bibitem{purbhoo}
K. Purbhoo, \textit{A Nullstellensatz For Amoebas}, preprint.

\bibitem{knots} W. H. Stevens, \emph{Recursion formulas for some abelian knot
invariants}, Journal of Knot Theory and Its Ramifications, Vol. 9,
No. 3 (2000) 413-422.

\bibitem{Zeil} D. Zeilberger, \emph{Dodgson's determinant-evaluation rule 
proved by two-timing men and women}, Elec. J. Comb. 4(2), 1997.

\end{thebibliography}
\end{document}